\documentclass[a4paper,11pt]{article}
\usepackage[mathscr]{eucal}
\usepackage{amsmath}
\usepackage{amssymb}
\usepackage{latexsym}
\usepackage{theorem}
\theoremstyle{plain}
\newtheorem{Theorem}{Theorem}[section]
\newtheorem{Corollary}[Theorem]{Corollary}

\newtheorem{Proposition}[Theorem]{Proposition}
\theorembodyfont{\rmfamily}
\newtheorem{Remark}[Theorem]{Remark}
\newtheorem{Example}[Theorem]{Example}
\newtheorem{Definition}[Theorem]{Definition}

\def\Vec#1{\mbox{\boldmath $#1$}}
\numberwithin{equation}{section}

\title{Deformations of surfaces preserving 
conformal or similarity invariants}

\author{Atsushi Fujioka and  Jun-ichi Inoguchi}
\date{\it Dedicated to professor Hideki Omori}

\begin{document}
\maketitle

\section*{Introduction}

In \cite{BPP}, Burstall, Pedit and Pinkall gave a 
fundamental theorem of surface theory in 
M\"obius 3-space in modern formulation.
Surfaces in M\"obius 3-space are determined by 
conformal Hopf differential and Schwartzian derivative
up to conformal transformations.
Isothermic surfaces are characterized as surfaces 
in M\"obius 3-space which admit
deformations preserving the conformal Hopf differential.

\vspace{0.2cm}

On every surface in M\"obius 3-space, a (possibly singular)
conformally invariant
Riemannian metric is introduced.  This metric is called the
{\it M\"obius metric} of the surface.
The Gaussian curvature of the M\"obius metric is
called the M\"obius curvature.
Here we point out that 
the preservation of M\"obius metric is weaker than that of conformal 
Hopf differential.

\vspace{0.2cm}

Constant mean curvature surfaces (abbriviated as CMC surfaces)
in the space forms are typical
examples of isothermic surfaces.
Bonnet showed that every constant mean 
curvature surface
admits 
a one-parameter family of isometric
deformations preserving the mean curvature.
A surface which admits such a family of
deformations is called a {\it Bonnet surface}.
Both the isothermic surfaces
and Bonnet surfaces  
are regarded as {\it geometric} generalizations of constant mean
curvature surfaces.

On the other hand, from the viewpoint of 
integrable system theory, Bobenko introduced the notion of
surface with harmonic inverse mean curvature
(HIMC surface, in short) in Euclidean 3-space $\mathbb{R}^3$. 
The first named author extended the notion of HIMC surface
in $\mathbb{R}^3$
to that of  3-dimensional space forms \cite{F}.
HIMC surfaces have deformation families (associated family)
which preserve the conformal structure of the surface and 
the harmonicity of the reciprocal mean curvature.
There exist local bijective conformal correspondences
between HIMC surfaces in different space forms.

It should be remarked that 
while Bonnet surfaces are isothermic,
HIMC surfaces are not necessarily isothermic.
In fact, the associated family
of a Bonnet surface or a HIMC surface preserves the 
M\"obius metric, while
the conformal Hopf differential
of a HIMC surface is not preserved 
in the associated family.
   
\vspace{0.2cm}

These observations motivate us to study surfaces in 
M\"obius 3-space (or space forms)
which admit deformations preserving the M\"obius metric.
We call such surfaces {\it M\"obius applicable surfaces}.

In this paper we study M\"obius applicable surfaces. 

First, we shall show the following new characterization of
Willmore surfaces (Theorem \ref{Willmore}):

\vspace{0.15cm}

\noindent
{\it A surface in M\"obius 3-space is Willmore if and only if it is 
a M\"obius applicable surface whose deformation familiy 
preserves the Schwarzian derivative}.

\vspace{0.2cm}

Next, we shall characterize both Bonnet surfaces and HIMC surfaces
in the class of M\"obius applicable surfaces 
in terms of similarity invariants
(Theorem \ref{BonnetHIMC}):

\vspace{0.15cm}

\noindent
{\it A surface in Euclidean 3-space is 
a Bonnet surface or a HIMC surface 
if and only if it is a M\"obius applicable surface
with specific deformation family in which,
the ratio of principal curvatures
is preserved}.

\vspace{0.15cm}

Furthermore we shall give the following characterization of
flat Bonnet surfaces (Theorem \ref{FlatBonnet}):

\vspace{0.15cm}

\noindent
{\it A Bonnet surface of 
non-constant mean curvature in Euclidean 3-space is 
flat if and only if its ratio of 
principal curvatures or 
M\"obius curvature is constant}.

\vspace{0.15cm}

Our characterization results imply that
``Bonnet" and ``HIMC" are similarity notion.
Thus these classes of surfaces fit naturally 
into similarity geometry.

We emphasize that similarity 
geometry provide
us 
non-trivial differential geometry of
integrable surfaces. In fact, 
the Burgers hierarchy are derived as
deformation of plane curves 
in similarity geometry.

\section{Deformation of surfaces preserving conformal invariants}

\subsection{Generalities of surface theory in conformal geometry}

Let $\mathbb{R}^3$ be the Euclidean 3-space.
The group $\mathrm{Conf}(3)$ of all
conformal diffeomorphisms are generated by 
isometries, dilations and inversions.
The conformal compactification $\mathcal{M}^{3}$ of 
$\mathbb{R}^3$ is called the {\it M\"obius $3$-space}.
By definition, $\mathcal{M}^3$ is the $3$-sphere equipped with
the canonical flat conformal structure.

In this paper, we use 
the projective lightcone model
of the M\"obius 3-space introduced by Darboux.

Let $\mathbb{R}^5_1$, be the {\it Minkowski $5$-space}
with canonical Lorentz scalar product:
$$
\langle \xi,\eta
\rangle=-\xi_{0}\eta_{0}+
\xi_{1}\eta_{1}+
\xi_{2}\eta_{2}+
\xi_{3}\eta_{3}+
\xi_{4}\eta_{4}.
$$
We denote the natural basis of $\mathbb{R}^5_1$ by
$\{\Vec{e}_{0},\Vec{e}_{1},\cdots,\Vec{e}_{4}\}$.
The unit timelike vecor $\Vec{e}_0$ time-orients $\mathbb{R}^5_1$.
The linear isometry group of $\mathbb{R}^{5}_{1}$ is denoted 
by $\mathrm{O}_{1}(5)$ and called the {\it Lorentz group}
\cite{O'Neill}.
The {\it lightcone} $\mathcal{L}$ of $\mathbb{R}^5_1$ is
$$
\mathcal{L}=\{v\in \mathbb{R}^{5}_{1} \
\vert \
\langle v,v 
\rangle=0,\
v\not=0\
\}.
$$
The lightcone has two connected components
$$
\mathcal{L}_{\pm}:=\{ v\in \mathcal{L} \ 
\vert \
\pm 
\langle \Vec{e}_{0},v\rangle<0\}.
$$
These connected components 
$\mathcal{L}_{+}$ and 
$\mathcal{L}_{-}$ are called 
the {\it future lightcone} and
{\it past lightcone},
respectively.

For $v\in \mathcal{L}$ and $r \in \mathbb{R}^{\times}$,
clearly, $rv \in \mathcal{L}$. Thus
$\mathbb{R}^{\times}$ acts freely on $\mathcal{L}$. The
quotient $\mathbb{P}(\mathcal{L})$ of $\mathcal{L}$ by 
the action of $\mathbb{R}^{\times}$ is
called the {\it projective lightcone}. 

The projective lightcone has a conformal structure
with respect to which it
is conformally equivalent
to the unit sphere $S^3$ with constant curvature $1$ metric. 

In fact, let us take a 
unit timelike vector $t_0$ and set
$$
S_{t_0}:=\{v\in \mathbb{P}(\mathcal{L})
\
\vert
\
\langle t_{0},v
\rangle=-1
\}.
$$
For $v \in S_{t_0}$, express $v$ as 
$v=v^{\perp}+t_0$ so that
$v^{\perp} \perp t_0$. Then 
$$
0=\langle v,v 
\rangle=\langle v^\perp,
v^\perp
\rangle+
\langle t_0,t_0
\rangle=
\langle v^\perp,
v^\perp
\rangle-1.
$$
This implies that 
the projection $v\mapsto v^\perp$ is an isometry
from $S_{t_0}\subset \mathbb{P}(\mathcal{L})$ onto the unit 3-sphere
$S^3$ in the Euclidean $4$-space $\mathbb{R}^4=
(\mathbb{R}t_{0})^{\perp}$.
This identification induces the following identification:
$$
\mathcal{M}^{3}\to \mathbb{P}(\mathcal{L}); \
v\longmapsto [1:v]
$$
between the M\"obius 3-space and the
projective lightcone.

More generally, all space forms are realized as
conic sections of $\mathcal{L}$. In fact, 
for a non-zero vector $v_0$, the section $S_{v_0}$
inherits a Riemannian metric of constant curvature 
$-\langle v_0,v_0\rangle$.

\begin{Definition}
A diffeomorphism of $\mathcal{M}^3$ is said to be a
{\it M\"obius transformation}
if it preserves $2$-spheres. The Lie group 
$\textrm{M\"ob}(3)$ of M\"obius transformations
is called the {\it M\"obius group}.
\end{Definition}

Any conformal diffeomorphism of $\mathcal{M}^3$
is a M\"obius transformation.

The following result is due to Liouville:
\begin{Proposition}
Let $\phi:U\to V$ be a conformal
diffeomorphism 
between two connected open subsets
of $\mathcal{M}^3$. Then there exists 
a unique M\"obius transformation $\Phi$ of $\mathcal{M}^3$
such that $\phi=\Phi\vert_{U}$.
\end{Proposition}

The linear action of Lorentz group
$\mathrm{O}_{1}(5)$
on $\mathbb{R}^{5}_{1}$
preserves $\mathcal{L}$ and decends to an action
on $\mathbb{P}(\mathcal{L})$.
Take a unit timelike vector $t_0$
and $T\in \mathrm{O}_{1}(5)$.
Then the restriction of $T$
to $S_{t_0}$ gives an isometry
 from $S_{t_0}$ onto $S_{Tt_0}$.
 This isometry $S_{t_0}\to S_{Tt_0}$
induces a conformal diffeomorphism 
on $\mathbb{P}(\mathcal{L})$.
These facts together with Liouville's theorem
imply that the following sequence
$$
0\to \mathbb{Z}_{2}\to \mathrm{O}_{1}(5)\to
\textrm{M\"ob}(3) 
\to 0
$$
is exact. Hence
$\textrm{M\"ob}(3)\cong \mathrm{O}^{+}_{1}(5)$,
where $\mathrm{O}_{1}^{+}(5)$ is the subgroup
of $\mathrm{O}_{1}(5)$ that preserves $\mathcal{L}_{\pm}$.
(See \cite[Theorem 1.2, 1.3]{Burstall})


The {\it de Sitter $4$-space}
$$
S^4_1=\{ v\in \mathbb{R}^{5}_{1} \ \vert \ \langle v,v \rangle=1\}
$$
parametrizes the space of all oriented 
conformal $2$-spheres in $\mathcal{M}^3$.
In fact, take a unit spacelike vector
$v\in S^4_1$ and denote by $V$ the 
$1$-dimensional linear subspace spanned by $v$. Then 
$\mathbb{P}(\mathcal{L}\cap V^{\perp})$ 
is a conformal $2$-sphere in $\mathcal{M}^3$.
Conversely  any conformal $2$-sphere can be represented in this form. 
Via this correspondence, the space of all conformal $2$-spheres
is identified with $S^4_1/\mathbb{Z}_{2}$.
Viewed as a surface $S_{v_0}\cap V^{\perp}$ of the conic section
$S_{v_0}$, this conformal $2$-sphere has the
mean curvature vector $\mathbb{H}_v$ 
$$
\mathbb{H}_{v}=-v_{0}^{\perp}-\langle v_0^\perp,
v_0^\perp\rangle v
$$
at $v$,
where $v_0$ is decomposed as 
$v_0=v_0^T+v_0^\perp$ according to the 
orthogonal direct sum $\mathbb{R}^5_1=V\oplus V^\perp$.

Let $F:M\to \mathcal{M}^3=
\mathbb{P}(\mathcal{L})$
be a conformal immersion of a 
Riemann surface into the M\"obius $3$-space.
The {\it central sphere congruence} (or {\it mean curvature
sphere}) of $F$ is a map $S:M\to S^4_1$ 
which assigns to each point $p\in M$, the
unique oriented $2$-sphere $S(p)$ tangent to 
$F$ at $F(p)$ which 
has the same orientation to $M$ and 
the same mean curvature vector
$\mathbb{H}_{S(p)}=\mathbb{H}_{p}$ at $F(p)$ as $F$.
The pull-back $\mathrm{I}_{\mathcal{M}}:=\langle dS,dS\rangle$
of the metric of $S^4_1$ by the central sphere congruence gives a
(possiblly singular) metric on $M$ and called the {\it M\"obius metric}
of $(M,F)$. The M\"obius metric is singular at umbilics.
The area functional $\mathcal{A}_{\mathcal{M}}$ of 
$(M,\mathrm{I}_{\mathcal{M}})$
is called the 
{\it M\"obius area} of $(M,F)$.
A conformally immersed surface $(M,F)$ is said to be a
{\it Willmore surface} if it is a critical point of the 
M\"obius area functional.

\subsection{The integrability condition}

Let $F:M \to \mathcal{M}^3$ be a 
conformal immersion of a Riemann surface.
A {\it lift} of $F$ is a map $\psi:M \to \mathcal{L}_{+}$ 
into the future lightcone such that
$\mathbb{R}\psi(p)=F(p)$ for any $p\in M$.
For instance, 
$\phi:=(1,F):M \to S_{\Vec{e}_0}\subset 
\mathcal{L}_{+}$ is a lift
of $F$.
This lift is called the {\it Euclidean lift} of $F$.
Now let $\phi$ be the Euclidean lift of $F$.
Then for any positive function $\mu$ on $M$,
$\phi \mu$ is still a lift of $F$. 
Direct computation shows that
$$
\langle d(\phi \mu),
d(\phi \mu)
\rangle_{1}=\mu^{2}
\langle dF,dF
\rangle_{1},
$$  
where $\langle\cdot,\cdot
\rangle_{1}$ is the constant curvature $1$ metric of $\mathcal{M}^3$.
Take a local complex coordinate $z$. Then the 
{\it normalized lift} $\psi$ with respect to $z$ is defined 
by the relation:
$$
\langle d\psi,d\psi
\rangle=dzd\bar{z}.
$$
This lift is M\"obius invariant. 
For another local complex coordinate $\tilde{z}$, the
normalized lift $\tilde{\psi}$ with respect to
$\tilde{z}$ is computed as $\tilde{\psi}=\psi|\tilde{z}_{z}|$. 

The normalized lift $\psi$ satisfies the following
inhomogeneous Hill equation: 
$$
\psi_{zz}+\frac{c}{2}\psi=\kappa.
$$
Under the coordinate change $z\mapsto \tilde{z}$,
the coefficients $c$ and $\kappa$ are changed as
$$
\widetilde{\kappa}\frac{d\tilde{z}^2}{|d\tilde{z}|}=
\kappa\frac{dz^2}{|dz|}.
$$
\begin{equation}\label{SchwarzianTransRule}
\tilde{c}d\tilde{z}^2=(c-S_{z}(\tilde{z}))dz^2,
\end{equation}
where $S_{z}(\tilde{z})$ is the {\it Schwarzian derivative}
of $\tilde{z}$ with respect to $z$.
Here we recall that the Schwarzian derivative
$S_{z}(f)$ of a meromorphic function
$f$ on $M$ is defined by
$$
S_{z}(f):=\left(
\frac{f_{zz}}{f_z}
\right)_{z}-\frac{1}{2}
\left(
\frac{f_{zz}}{f_{z}}
\right)^{2}.
$$ 
Moreover two meromorphic functions $f$ and $g$ are
{\it M\"obius equivalent}, i.e.,
related by a linear fractional transformation:
$$
g=\frac{af+b}{cf+d},
\ \ 
 \left(\begin{array}{cc}
a& b\\
c & d
\end{array}
\right)\in
\mathrm{SL}_{2}\mathbb{C}
$$
if and only if their Schwartzian derivatives $S_{z}(f)=S_{z}(g)$
agree.

Now we denote by $L$ the {\it $1$-density bundle} of $M$:
$$
L:=(K \otimes_{\mathbb C}\overline{K})^{-\frac{1}{2}},\ 
K \ \textrm{is the canonical bundle of } M. 
$$
The transformation law (\ref{SchwarzianTransRule}) 
implies that $\kappa \> dz^2/|dz|$
is a section of $L\>K^2$, i.e., an $L$-valued quadratic
differential on $M$. This section is called the {\it conformal Hopf differential} of $(M,F)$. The differential $cdz^2$ is called the
{\it Schwarzian} of $(M,F)$. 
The coefficient function $c$ is also called the Schwarzian.
 
Note that the conformal Hopf differential vanishes
identically if and only if $M$ is
totally umbilical.

The integrability condition for a conformal
immersion $F:M\to \mathcal{M}^{3}$ is
given in terms of $\kappa$ and $c$ as follows:
\begin{equation}
\begin{cases}
\frac{1}{2}
c_{\bar{z}}=
3\bar{\kappa}_{z}
\kappa+\bar{\kappa}\kappa_z,
\\
\textrm{Im}\hspace{0.5mm}
\left(\kappa_{\bar{z}\bar{z}}+
\frac{1}{2}\bar{c}\kappa
\right)=0.
\end{cases}
\end{equation}
These equations are called, the 
{\it conformal Gauss equation}
and the {\it conformal Codazzi equation},
respectively.

The M\"obius metric $\mathrm{I}_{\mathcal{M}}$ is 
represented by
\begin{equation}
\mathrm{I}_{\mathcal{M}}=4|\kappa|^{2}dzd{\bar z}.
\end{equation}

The Euler-Lagrange equation for the 
M\"obius area functional $\mathcal{A}_{\mathcal{M}}$ 
is called the {\it Willmore surface equation} and given 
in terms of $\kappa$ and $c$ as follows 
if \cite[p.~51]{BPP}:
 
\begin{equation}
\kappa_{\bar{z}\bar{z}}+
\frac{1}{2}\bar{c}\kappa
=0.
\end{equation}

\subsection{Deformation of surfaces preserving the Schwarzian \\
derivative or the conformal Hopf differential}

Generally speaking, the conformal Hopf
differential alone determines surfaces in $\mathcal{M}^3$.
However, there are the only exceptional surfaces--{\it isothermic surfaces}
\cite{BC}.
Isothermic surfaces are defined as surfaces in $\mathcal{M}^{3}$
conformaly parametrized by their curvature lines away from
umbilics.
Away from umbilics, there are
holomorphic coordinates in which
the conformal Hopf
differential is real valued.
Such holomorphic coordinates (and their associated real
coordinates) are called {\it isothermic coordinates}.

Now let $(M,F)$ be an isothermic surface
parametrized by an isothermic coordinate $z$.
Then
under the deformation: 
$$
c\longrightarrow c_{r}:=c+r,\ \
r\in \mathbb{R},
$$
the conformal Gauss-Codazzi equations
$$
c_{\bar z}=4(\kappa^{2})_z,
\ \
\mathrm{Im}\> (\kappa_{\bar{z}\bar{z}}+\frac{1}{2}\bar{c}\kappa)=0
$$ 
are invariant.
Hence, as in the case of CMC surfaces,
one obtains a 1-parameter family $\{F_r\}$ of deformations
through $F=F_0$ preserving the conformal Hopf differential $\kappa$.
Since all $c_r$ are distinct,
the surfaces $\{F_r\}$ are non-congruent, each other.
The family $\{F_r\}$ is refered as the {\it associated family}
of an isothermic surface $(M,F)$. 
The correspondence $F\longrightarrow F_r$ 
is called the {\it $T$-transformation} by Bianchi
\cite{Bianchi1905b}.
The $T$-transformation was 
independently introduced by 
Calapso \cite{Cal5} and is also called 
the {\it Calapso transformation}. 

The existence of deformations preserving the
conformal Hopf differential characterizes 
isothermic surfaces as follows:

\begin{Theorem}{\rm(\cite{BPP})}
A surface in $\mathcal M^3$ is isothermic if and only if it has 
deformations 
preserving the conformal Hopf differential.
\end{Theorem}
\noindent

\begin{Corollary}{\rm(\cite{BPP})}
Let $F_{1},F_{2}:M\to \mathcal{M}^{3}$
be two non-congruent
surfaces with same conformal Hopf differential. Then
both $F_1$ and $F_2$ belong to the same associated family of
an isothermic surface.
\end{Corollary}

On the other hand, for deformations preserving M\"obius metric and
\newline
\noindent
Schwarzian, we have the
following {\it new} characterization of
Willmore surfaces.

\begin{Theorem}\label{Willmore}
A surface in $\mathcal M^3$ is Willmore if and only if it has 
M\"obius-isometric 
deformations preserving the Schwarzian derivative.
\end{Theorem}
\noindent
\textit{Proof.} Let $F$ be a surface in $\mathcal M^3$ with the Schwarzian derivative $c$ and 
the conformal Hopf differential $\kappa$. 
If $F$ has deformation preserving
the M\"obius metric $\mathrm{I}_{\mathcal{M}}$
and $c$, there exists an $S^1$-valued function $\lambda$ such that $\lambda\kappa$ and $c$ 
satisfy the conformal Gauss equation. 
Combining this with the conformal Gauss equation for $F$, we have 
\begin{equation}
3\bar{\lambda}_z\lambda+\bar{\lambda}\lambda_z=0,
\end{equation}
which implies that $\lambda^3\bar{\lambda}$ is holomorphic and hence $\lambda$ is an $S^1$-valued constant. 
Since $\lambda\kappa$ and $c$ satisfy the conformal Codazzi equation, combinig this with the conformal Codazzi equation 
for $F$, we have 
\begin{equation}
\kappa_{\bar{z}\bar{z}}+\frac{1}{2}\bar{c}\kappa=0,
\end{equation}
which implies that $F$ is Willmore. $\square$

\begin{Remark}
{\'E}.~Cartan formulated a general theory of
deformation of submanifolds in homogeneous
spaces. The classical deformation problems
(also called {\it applicability} of submanifolds
in classical literatures) in Euclidean,
projective and conformal geometry are
covered by Cartan's framework
\cite{Cartan1}-\cite{Cartan2}. 
 
According to Griffiths \cite{G} and Jensen \cite{Jensen},
two immersions 
$F_{1},\>F_{2}:M\to G/K$
of a manifold into a homogeneous manifold
are said to be {\it $k$-th order deformation} of each
other if there exists a smooth map $g:M\to G$ such that,
for every $p\in M$,
the 
Taylor expansions about 
$p$ of $F_2$ and $g(p)F_{1}$ agree through $k$-th order terms. 
An immersion $F:M\to G/K$ is said to be
{\it deformable of order $k$} if it
admits a non-trivial
$k$-th order deformation. 

Musso \cite{Musso} showed that a conformal immersion of a Riemann 
surface $M$ into the M\"obius $3$-space is $2$nd order deformable 
if and only if it is isothermic.
\end{Remark}

\begin{Remark}(Special isothermic surfaces)
Among isothermic surfaces in $\mathbb{R}^{3}$,
Darboux \cite{Darboux1899}
distinguished the class of 
{\it special isothermic surfaces}. 
An isothermic surface  $F:M\to \mathbb{R}^3$
with first and second fundamental forms;
$$
\mathrm{I}=e^{\omega}(dx^2+dy^2),
\ 
\mathrm{I}\!\mathrm{I}=e^{\omega}(k_{1}dx^2+k_{2}dy^{2})
$$
is called {\it special} of type $(A,B,C,D)$
if its mean curvature $H$ satisfies
the equation:
$$
4e^{\omega}|\nabla H|^{2}+m^{2}+2Am+2BH+2C\ell+D=0,
$$
where $\ell=2e^{\omega}\sqrt{H^2-K}$,
$m=-H\ell$ and $A,B,C,D$ are real constants.
Constant mean curvature surfaces are paticular examples of special
isothermic surface. 
Special isothermic surfaces with $B=0$  
are conformally invariant.
Moreover, Bianchi \cite{Bianchi1905} and 
Calapso \cite{Cal5} showed that
an umbilic free isothermic surface
in $\mathcal{M}^3$  is special with $B=0$ 
if and only if it is conformally equivalent to
a constant mean curvature surface in space forms.
For modern treatment of special isothermic surfaces and their Darboux
transformations, we refer to \cite{MN1}.
In \cite{Bernstein}, Bernstein constructed 
non-special, non-canal isothermic tori in
$\mathcal{M}^3$ with spherical curvature lines.
\end{Remark}

Let $F:M\to \mathcal{M}^3$ be a conformal immersion.
Then $F$ is said to be a {\it constrained Willmore surface}
if it is a critical point of the M\"obius area functional
under (compactly supported) conformal variations.

\begin{Proposition}{\rm(\cite{BPP})}
$F:M\to \mathcal{M}^{3}$ is constrained Willmore if and only
if there exists a holomorphic quadratic differential
$qdz^2$ such that
\begin{equation}\label{cWillmore}
\kappa_{\bar{z}\bar{z}}+\frac{1}{2}\bar{c}
\kappa=\mathrm{Re}\>
(\bar{q}\kappa).
\end{equation}
\end{Proposition}

The constrained Willmore surface
equation (\ref{cWillmore}) 
has the following
deformation:
$$
\kappa\longrightarrow \kappa_{\lambda}:=\lambda \kappa,
\
c\longrightarrow c_{\lambda}:=c+(\lambda^{2}-1)q,
\
q\longrightarrow q_{\lambda}:=
\lambda q,
$$
for $\lambda \in S^1$.

Hence we obtain a one-parametric conformal deformation family
$\{F_{\lambda}\}$ of a constrained Willmore surface $(M,F)$.
This familiy is refered as the {\it associated
family} of $F$.

Obviously, for Willmore surfaces $(q=0)$,
the associated family preserves the Schwarzian.

The following characterization of 
constrained Willmore surfaces 
can be verified in a way similar 
to the proof of Theorem \ref{Willmore}:
\begin{Proposition}
A surface $F:M\to \mathcal{M}^{3}$ has a deformation
of the form 
$$
\kappa\longrightarrow \lambda \kappa,
\
c\longrightarrow c+r
$$
for some $S^1$-valued function $\lambda$ and
a holomorphic quadratic differential
$rdz^2$
if and only if $M$ is a constrained Willmore surface.
\end{Proposition}

\begin{Remark}
A classical result by Thomsen says
that a surface is isothermic Willmore 
if and only if it is minimal in 
a space form (\cite{Blaschke}, \cite[Theorem 3.6.7]{Udo},
\cite{Thomsen}).
Constant mean curvature surfaces
in space forms are isothermic and
constrained Willmore.
Richter \cite{Richter} showed that  
in the case of immersed tori in $\mathcal{M}^{3}$,
every isothermic constrained Willmore tori are
constant mean curvature tori in some space forms.
In contrast to Thomsen's result,
the assumption ``tori" is essential for Richter's result.
In fact, Burstall construced isothermic constrained 
Willmore cylinders which are not realized as constant mean 
curvature surfaces in any space forms. See \cite{BoPePi}.
\end{Remark}


\section{Deformation of surfaces preserving similarity \\
invariants}

As we saw in the preceding section, preservation of  
conformal Hopf differential is a strong restriction in 
the study of deformation of surfaces. 
Clearly, the preservation of M\"obius metric is weaker than that of 
conformal Hopf differential. In this section we study 
deformation of surfaces preserving the M\"obius metric.  

\subsection{M\"obius invariants via metrical language}
First, we discuss relations between metrical invariants and M\"obius  
invariants.

Let $F:M\to \mathbb{R}^{3}$ be a conformal immersion 
of a Riemann surface into the Euclidean $3$-space. 
Denote by $\mathrm{I}$ the first fundamental form (induced metric)
of $M$.
The Levi-Civita connections $D$ of $\mathbb{R}^3$
and $\nabla$ of $M$ are
related by the {\it Gauss equation}:
$$
D_{X}F_{*}Y=F_{*}(\nabla_{X}Y)+\mathrm{I}\!\mathrm{I}(X,Y)\Vec{n}.
$$
Here $\Vec{n}$ is the unit normal vector field.
The symmetric tensor field 
$\mathrm{I}\!\mathrm{I}$ is the {\it second
fundamental form} derived from $\Vec{n}$.

The trace free part of the second fundamental form is
given by $\mathrm{I}\!\mathrm{I}-H\mathrm{I}$, where $H$ is the mean curvature
function. Define a function $h$ by $h:=\sqrt{H^2-K}$.
This function $h$ is called the 
{\it Calapso potential}. 

Then one can check that the normal vector field 
$\Vec{n}/h$ and the symmetric tensor field 
$h^2 \>\mathrm{I}$ are invariant under the conformal change of the
ambient Euclidean metric.
Moreover the trace free symmetric tensor field 
$$
\mathrm{I}\!\mathrm{I}_{\mathcal{M}}:=
h(\mathrm{I}\!\mathrm{I}-H\mathrm{I})
$$
is also conformally invariant.
It is easy to see that $h^{2}\> \mathrm{I}$
coincides with the  
{\it M\"obius metric} $\mathrm{I}_{\mathcal{M}}$
of $(M,F)$. 
The pair 
$(\mathrm{I}_{\mathcal{M}},\mathrm{I}\!\mathrm{I}_{\mathcal{M}})$
is called the {\it Fubini's conformally invariant fundamental forms}.
The Gaussian curvature $K_{\mathcal{M}}$ of $(M,\mathrm{I}_{\mathcal{M}})$
is called the {\it M\"obius curvature} of $(M,F)$.
The M\"obius area functional $\mathcal{A}_{\mathcal{M}}$
of $(M,\mathrm{I}_{\mathcal{M}})$
is computed as
$$
\mathcal{A}_{\mathcal{M}}=\int_{M}(H^2-K)dA_{\mathrm{I}}.
$$

Now let us take a local complex coordinate $z$ and 
express the first fundamental form as 
$\mathrm{I}=e^{\omega}dzd\bar{z}$. The (metrical) 
{\it Hopf differential} is defined by
$$
Q^{\#}:=Qdz^{2},\ 
Q=\langle F_{zz},\Vec{n}\rangle.
$$
Then the conformal Hopf differential and the mestrical one are related
by the formula:
\begin{equation}\label{HopfdiffsRelation}
\kappa=Qe^{-\frac{\omega}{2}}.
\end{equation}

The Schwarzian derivative is
repressented as
$$
c=\omega_{zz}-\frac{1}{2}(\omega_{z})^{2}+2HQ.
$$

\subsection{Similarity geometry}

The similarity geometry is a subgeometry of M\"obius geometry
whose symmetry group is the {\it similarity transformation group}:
$$
\mathrm{Sim}(3)=\mathrm{CO}(3)\ltimes \mathbb{R}^{3},
$$
where $\mathrm{CO}(3)$ is the linear conformal group
$$
\mathrm{CO}(3)=\{A \in \mathrm{GL}_{3}\mathbb{R}
\
\vert
\
{}^{\exists} c \in \mathbb{R};
\
{}^{t}AA=cE \ \}.
$$

Let $F:M\to \mathbb{R}^3$ be an immersed surface with 
unit normal $\Vec{n}$ as before.

Under the similarity transformation of $\mathbb{R}^3$,
Levi-Civita connections 
$D$ and $\nabla$ are invariant. Hence the vector valued
second fundamental form $\mathrm{I}\!\mathrm{I}\Vec{n}$ is
similarity invariant.  
The shape opeartor $S=-d\Vec{n}$ itself is not similarity invariant,
but the ratio of principal curvatures are invariant.
It is easy to see that the constancy of
the ratio of principal curvatures is equivalent 
to the constancy of $K/H^2$. The function $K/H^2$ 
is similarity invariant.
The principal directions are yet another similarity invariant.

\subsection{Deformation of surfaces preserving 
the M\"obius metric and the ratio of principal curvatures}

Let $F:M\to \mathbb{R}^3$ be a surface in Euclidean $3$-space.
Then the Gauss-Codazzi equations of $(M,F)$ are given by
\begin{equation}\label{2.1.X}
\begin{cases}
\omega_{z\bar{z}}+\frac{1}{2}H^2e^{\omega}-2\vert Q\vert^2e^{-\omega}=0,\\
Q_{\bar{z}}=\frac{1}{2}H_ze^{\omega}.
\end{cases}
\end{equation}

The Gauss-Codazzi equations imply the following fundamental fact
due to Bonnet.

\begin{Proposition}{\rm(\cite{Bonnet})}
Every non-totally umbilical constant mean curvature surface 
admits a one-parameter isometric deformation  
preserving the mean curvature.
\end{Proposition}

Here we exhibit two examples of surfaces
which admit deformations preseving the M\"obius
metric.

\begin{Example}{\rm(Bonnet surfaces)}
Let $F:M\to \mathbb{R}^3$ be a Bonnet 
surface. Namely $(M,F)$ admits  
a non-trivial isometric deformation
$F\longrightarrow F_\lambda$
preserving the mean curvature. 
The deformation family 
$\{F_\lambda\}$ is called the
{\it associated family} of $(M,F)$.

Since  all the members $F_\lambda$ have the same metric and the mean
curvature, they have the same M\"obius metric.  
Note that the conformal Hopf differential is not 
preserved under the deformation.
\end{Example}

\begin{Example}{\rm(HIMC surfaces)}
A surface $F:M\to \mathbb{R}^3$
is said 
to be a {\it surface with harmonic inverse mean curvature}
({\it HIMC surface}, in short)
if its inverse mean curvature function $1/H$
is a harmonic function on $M$ \cite{Bo}.
Since $1/H$ is harmonic, $H$ can be expressed as 
$1/H=h+\bar{h}$ for some holomorphic function $h$.
The associated family
$\{F_\lambda\}$ of a HIMC surface $F$ is given by
the following
metrical data
$(\mathrm{I}_{\lambda},H_{\lambda},Q_{\lambda})$:
$$
\mathrm{I}_{\lambda}=e^{\omega_\lambda}dzd\bar{z},
\ \
e^{\omega_\lambda}=\frac{e^\omega}
{(1-2\sqrt{-1}\>\bar{h}t)^{2}(1+2\sqrt{-1}ht)^{2}},
$$
$$
\frac{1}{H_\lambda}=h_{\lambda}+\overline{h_\lambda},
\ \
h_{\lambda}=
\frac{h}{1+2\sqrt{-1}ht},
$$
$$
Q_{\lambda}=\frac{Q}{(1+2\sqrt{-1}ht)^{2}},
\ \
\lambda=
\frac{1-2\sqrt{-1}\bar{h}t}{1+2\sqrt{-1}ht},
\ \ t\in \mathbb{R}.
$$
From these,
we have
$$
(H_{\lambda}^{2}-K_{\lambda})=(1-2\sqrt{-1}\bar{h}t)(1+2\sqrt{-1}ht)(H^2-K),
$$
Hence
$$
(H_{\lambda}^{2}-K_{\lambda})e^{\omega_\lambda}=
(H^2-K)e^{\omega}.
$$
Thus the M\"obius metric is preserved under the
deformation $F\longrightarrow F_\lambda$.
On the other hand, 
the conformal Hopf differential
is not preserved under the deformation.
In fact, the conformal Hopf
differential 
of $F_\lambda$ is
$$
\kappa_{\lambda}:=Q_{\lambda}e^{-\frac{\omega_\lambda}{2}}
=\kappa\>
\frac{1-2\sqrt{-1}\bar{h}t}{1+2\sqrt{-1}ht}.
$$
Clearly $\kappa_\lambda$ is not preserved
under the deformation.
%

While Bonnet surfaces are isothermic, HIMC surfaces
are not necesarrily so.
The dual surfaces of Bonnet surfaces are isothermic HIMC
surfaces. 
Since the associated families of Bonnets surface or 
isothermic HIMC surfaces do not preserve the conformal
Hopf differential, these families differ from the $T$-transformation
families. Note that $T$-transformations are only well defined up to
M\"obius transformations \cite[section 2.2.3]{Burstall}.
\end{Example}

\vspace{0.2cm}

Now we prove the following theorem which characterizes Bonnet surfaces and
HIMC surfaces in the class of surfaces which posses 
M\"obius metric preserving deformations.
We call such surfaces {\it M\"obius applicable surfaces}.

\begin{Theorem}\label{BonnetHIMC}  
Let $F$ be a surface in 
$\mathbb{R}^3$ which has deformation preserving the 
M\"obius metric and the ratio 
of principal curvatures. 
Then the deformation is given by 
\begin{equation}
e^{\omega}\to\vert\lambda\vert^2 e^{\omega},\ H\to\frac{1}{\vert\lambda\vert}H,\ Q\to\lambda Q,
\end{equation}
where $\lambda$ is a function with $\vert\lambda\vert=\vert f\vert$ for some holomorphic function $f$. 
Moreover if $\vert\lambda\vert=1$ (respectively $\lambda$ is holomorphic), then $F$ is a Bonnet surface 
(respectively a HIMC surface). 
\end{Theorem}
\noindent
\textit{Proof.} 
Note that the quantities $\vert Q\vert^2e^{-\omega}$ and 
$e^{-\omega}/H^2$ are 
invariant under the deformation, 
which implies that the deformation is given by 
as above for some function $\lambda$
(see (\ref{HopfdiffsRelation})). 
From the Gauss equation we have 
\begin{equation*}
(\log\vert\lambda\vert^2)_{z\bar{z}}=0,
\end{equation*}
which implies that $\vert\lambda\vert=\vert f\vert$ for some holomorphic function $f$. 

If $\vert\lambda\vert=1$, the deformation is nothing but the isometric deformation preserving the mean curvature. 
Hence $F$ is a Bonnet surface. 

If $\lambda$ is holomorphic, putting $(H')^2=H^2/\vert\lambda\vert^2$ and differentiating it by $z$, we have 
\begin{equation*}
2H'H'_z=-\frac{\bar{\lambda}\lambda_z}{\vert\lambda\vert^4}H^2+\frac{2}{\vert\lambda\vert^2}HH_z.
\end{equation*}
Note that $Q\not=0$ since $F$ is umbilic-free. 
Combining the Codazzi equations for $F$ and the surface obtained by deformation, we have 
\begin{equation}\label{2.1.y}
H'_z=\frac{1}{\bar{\lambda}}H_z.
\end{equation}
Hence we have
\begin{equation*}
H'=-\frac{\lambda_z}{2\lambda^2}\frac{H^2}{H_z}+\frac{1}{\lambda}H.
\end{equation*}
Differentiating it by $\bar{z}$ and using 
(\ref{2.1.y}) again, we have
\begin{equation*}
H_{z\bar{z}}-\frac{2\vert H_z\vert^2}{H}=0,
\end{equation*}
which implies that $F$ is a HIMC surface. $\square$

\subsection{Flat Bonnet surfaces}

Let $M$ be a Bonnet surface in $\mathbb{R}^3$.
Then away from umbilics, there exists an 
isothermic coordinate $z$ such that the Gauss-Codazzi
equations of $M$ reduces to the following
third order  
ordinary differential equation ({\it Hazzidakis equation} \cite{Hazzidakis}):

\begin{equation}\label{2.2.x}
\left\{\left(\frac{H_{ss}}{H_s}\right)_s-H_s\right\}R^2=2-\frac{H^2}{H_s},
\ \ \
H_{s}<0,
\end{equation}
where $s=z+\bar{z}$
and the coefficient function $R(s)$ is 
one of the following function \cite[p.~30]{BE}:
$$
R_{A}(s)=\frac{\sin (2s)}{2},\ \
R_{B}(s)=\frac{\sinh (2s)}{2},
\ \
R_{C}(s)=s.
$$
The modulas $|Q|$ of the metrical Hopf
differentai $Qdz^2$ is given by
\begin{equation}\label{Hopfdiffs-Bonnet}
|Q(z,\bar{z})|=\frac{1}{R(s)^{2}}.
\end{equation} 
 A Bonnet surface is said to be 
 of {\it type} $A$, $B$ or $C$,
 respectively if away from
 critical points of the mean curvature,
 it is determined by 
a solution to  Hazzidakis equation with coefficient
$R_{A}$, $R_{B}$ or $R_{C}$
(\cite[Definition 3.2.1]{BE}, \cite{Cartan3}).

\begin{Proposition}{\rm(\cite{BE}, \cite{FI})}
Flat Bonnet surfaces in $\mathbb{R}^3$ are
of $C$-type.
\end{Proposition}

Flat Bonnet surfaces are characterized as follows 
in terms of conformal (M\"obius) or similarity invariants.

\begin{Theorem}\label{FlatBonnet} 
A Bonnet surface in $\mathbb{R}^3$ 
with non-constant mean curvature 
is flat if and only if the
M\"obius curvature or
the ratio of 
the principal curvatures is constant. 
\end{Theorem}
\medskip

\noindent
\textit{Proof.} 
First we consider Bonnet surfaces with constant ratio 
of principal curvatures.
By the assumption the function $K/H^2$ is constant.
Computing $K/H^2$ by using (\ref{2.2.x})
and (\ref{Hopfdiffs-Bonnet}),
one can deduce that $K=0$ if $K/H^2$ is constant.

\vspace{0.2cm}

Next, 
the M\"obius curvature 
$K_{\mathcal M}$ is computed as  
\begin{equation*}
K_{\mathcal M}=\frac{1}{H_s}(\log H_s)_{ss}
\end{equation*}
by using the Hazzidakis equation (\ref{2.2.x}).

If $K_{\mathcal M}$ is constant, 
a direct computation shows that the 
solution of (\ref{2.2.x}) is 
\begin{equation*}
H=-\frac{2}{K_{\mathcal M}}\frac{1}{s}
\end{equation*}
with $K_{\mathcal M}<0$. Hence the surface is flat. $\square$

\section*{Appendices}
\subsection*{ A.1.
Curves in similarity geometry.}
%

Let us consider plane curve geometry in the
$2$-dimensional similarity geometry
$(\mathbb{R}^2,\mathrm{Sim}(2))$.
Here $\mathrm{Sim}(2)$ denotes the similarity
transformation group of $\mathbb{R}^2$.

Let $\gamma(s)$  be a regular curve on $\mathbb{R}^2$
parametrized by the Euclidean arclength $\sigma$.
Then the $\mathrm{Sim}(2)$-invariant parameter $s$ is
the {\it angle function}
$s=\int^\sigma \kappa_{E}(\sigma)d\sigma$, where
$\kappa_{E}$ is the Euclidean curvature function.
The $\mathrm{Sim}(2)$-invariant curvature
$\kappa_{S}$ is given by
$\kappa_{S}=(\kappa_{E})_{\sigma}/\kappa_E^2$.
Obviously, every circle is a curve of similarity curvature $0$.
The $\mathrm{Sim}(2)$-invariant frame field 
$\mathcal{F}=(T,N)$ is given by
$$
T=\gamma_{s},\ N=T_{s}+\kappa_{S} T.
$$
The Frenet-Serret equation of $\mathcal{F}$ is
$$
\mathcal{F}^{-1}\frac{d\mathcal{F}}{ds}=\left(
\begin{array}{cc}
-\kappa_{S} & -1\\
1 & -\kappa_{S}
\end{array}
\right).
$$

Now let us consider plane curves 
of nonzero constant similarity curvature. 

Put
$\kappa_{S}=c_1 \>(\textrm{constant})$.
Then we have
$1/\kappa_E=(-c_1)\sigma+c_{2}$, 
namely 
$\gamma$ is a curve whose inverse Euclidean curvature
$1/\kappa_E$
is a linear function of the Euclidean arclength parameter. 
Thus $\gamma$ is a log-spiral (if $c_1\not=0$)
or a circle ($c_1=0, c_2\not=0$).

These curves provide fundamental
examples of Bonnet surfaces as well as
HIMC surfaces.
In fact, let $\gamma$ be a plane curve of constant similarity curvature.
Then cylinder over $\gamma$ is a flat Bonnet 
surface in $\mathbb{R}^3$
as well as a flat HIMC surface in $\mathbb{R}^3$.
Generally, the Hazzidakis equation of Bonnet or isothermic HIMC 
surfaces reduces to Painlev{\'e} equations of type III, V or VI.
The solutions to log-sprial cylinder are elementary function solutions
to these Painlev{\'e} equations. 
(see \cite{BE}, \cite{FI}).

\subsection*{A.2. Time evolutions}

Let us consider the time evolution of a plane curve 
$\gamma(s)$ in similarity geometry.

Denote by $\gamma(s;t)$ the time evolution
which preserves the similarity arclength parameter $s$;
 
$$
\frac{\partial}{\partial t}\>\gamma(s;t)=gN+fT
$$
Then the similarity curvature $u=\kappa_{S}$ obeys the following
partial differential equation:
$$
u_{t}=f_{sss}-2uf_{ss}
-(3u_{s}-u^{2}-1)f_{s}
-(u_{ss}-2uu_{s})f+au_{s},
\ a\in \mathbb{R}.
$$
In particular, if we choose 
$f=-1,a=0$, then the
time evolution of $u$ 
obeys the 
{\it Burgers equation}:
$$
u_{t}=u_{ss}-2uu_{s}.
$$
More generally, the Burgers hierarchy is induced by
the above time evolution, see \cite[pp.~17--18]{CQ1}. 
Space curves in similarity geometry and their time evolution, we 
refer to \cite{CQ2}.

\vspace{0.5cm}

(A.~F.) Graduate School of Economics,
Hitotsubashi University, 2-1,
Naka, Kunitachi, Tokyo,
186--8601,
Japan

{\tt fujioka@math.hit-u.ac.jp}

\vspace{0.2cm}

(J.~I.) Department of Mathematics Education,
Utsunomiya University, 
Utsunomiya, 321-8505,
Japan

{\tt inoguchi@cc.ustunomiya-u.ac.jp}

\end{document}